\documentclass[a4paper,12pt]{article}
\usepackage{a4wide}
\usepackage{amsmath}
\usepackage{amssymb}
\usepackage{amsthm}
\usepackage{latexsym}
\usepackage{graphicx}
\usepackage[english]{babel}
\usepackage{makeidx}

\newtheorem{obs} [subsection]{Remark}
\newtheorem{exm} [subsection]{Example}

\newtheorem{prop}[subsection]{Proposition}

\newtheorem{teor}[subsection]{Theorem}
\newtheorem{lema}[subsection]{Lemma}
\newtheorem{cor} [subsection]{Corollary}
\newcommand{\Zng}{$\mathbb Z^n$-graded $S$-module}

\def\sdepth{\operatorname{sdepth}}
\def\depth{\operatorname{depth}}
\begin{document}
\selectlanguage{english}
\frenchspacing

\large
\begin{center}
\textbf{Stanley depth of monomial ideals with small number of generators}

Mircea Cimpoea\c s
\end{center}
\normalsize

\begin{abstract}
For a monomial ideal $I\subset S=K[x_1,\ldots,x_n]$, we show that $\sdepth(S/I)\geq n-g(I)$, where $g(I)$ is the number of the minimal monomial generators of $I$. If $I=vI'$, where $v\in S$ is a monomial, then we see that $\sdepth(S/I)=\sdepth(S/I')$. We prove that if $I$ is a monomial ideal $I\subset S$ minimally generated by three monomials, then $I$ and $S/I$ satisfy the Stanley conjecture. Given a saturated monomial ideal $I\subset K[x_1,x_2,x_3]$ we show that $\sdepth(I)=2$. As a consequence, $\sdepth(I)\geq \sdepth(K[x_1,x_2,x_3]/I)+1$ for any monomial ideal in $I\subset K[x_1,x_2,x_3]$.


\end{abstract}

\section*{Introduction}

Let $K$ be a field and $S=K[x_1,\ldots,x_n]$ the polynomial ring over $K$.
% We denote $\me=(x_1,\ldots,x_n)$ the maximal graded ideal of $S$.
Let $M$ be a \Zng. A \emph{Stanley decomposition} of $M$ is a direct sum $\mathcal D: M = \bigoplus_{i=1}^rm_i K[Z_i]$ as $K$-vector space, where $m_i\in M$, $Z_i\subset\{x_1,\ldots,x_n\}$ such that $m_i K[Z_i]$ is a free $K[Z_i]$-module. We define $\sdepth(\mathcal D)=min_{i=1}^r |Z_i|$ and $\sdepth(M)=max\{\sdepth(M)|\;\mathcal D$ is a Stanley decomposition of $M\}$. The number $\sdepth(M)$ is called the \emph{Stanley depth} of $M$. Herzog, Vladoiu and Zheng show in \cite{hvz} that this invariant can be computed in a finite number of steps if $M=I/J$, where $J\subset I\subset S$ are monomial ideals.

There are two important particular cases. If $I\subset S$ is a monomial ideal, we are interested in computing $\sdepth(S/I)$ and $\sdepth(I)$. There are some papers regarding this problem, like \cite{hvz},\cite{asia},\cite{sum}, \cite{shen} and \cite{mir}. Stanley's conjecture says that $\sdepth(S/I)\geq \depth(S/I)$, or in the general case, $\sdepth(M)\geq \depth(M)$, where $M$ is a finitely generated multigraded $S$-module. The Stanley conjecture for $S/I$ was proved for $n\leq 5$ and in other special cases, but it remains open in the general case. See for instance, \cite{apel}, \cite{hsy}, \cite{jah}, \cite{pops}, \cite{popi} and \cite{pope}.

Let $I\subset S$ be a monomial ideal. We assume that $G(I)=\{v_1,\ldots,v_m\}$, where $G(I)$ is the set of minimal monomial generators of $I$. We denote $g(I)=|G(I)|$, the number of minimal generators of $I$. Let $v=GCD(u|\; u\in G(I))$. It follows that $I=vI'$, where $I'=(I:v)$. For a monomial $u\in S$, we denote $supp(u)=\{x_i:\;x_i|u\}$. We denote $supp(I)=\{x_i:\;x_i|u$ for some $u\in G(I)\}$. We denote $c(I)=|supp(I')|$. 
In the first section, we prove results regarding some relations between $\sdepth(S/I)$, $\sdepth(I)$, $g(I)$ and $c(I)$. 

%For a monomial ideal $I\subset S$, $I^{sat}=(I:(x_1,\ldots,x_n)^{\infty}) = \bigcup_{k\geq 1} (I:(x_1,\ldots,x_n)^k)$ is called the saturation of $I$. We prove that  $\sdepth(S/I)=0$ if and only if $I\neq I^{sat}$, see Theorem $1.5$. 
In the second section, we give some applications. We prove that a monomial ideal $I\subset S$ minimally generated by three monomials has $\sdepth(I)=n-1$, see Theorem $2.4$. We prove that if $I$ is a monomial ideal $I\subset S$ minimally generated by three monomials, then $I$ and $S/I$ satisfy the Stanley conjecture, see Theorem $2.6$. 

\footnotetext[1]{This paper was supported by CNCSIS, ID-PCE, 51/2007}

In the third section, we prove that if $I\subset K[x_1,x_2,x_3]$ is saturated, then $\sdepth(I)\geq 2$, see Proposition $3.1$. As a consequence, $\sdepth(I)\geq \sdepth(K[x_1,x_2,x_3]/I)+1$ for any monomial ideal in $I\subset K[x_1,x_2,x_3]$, see Corollary $3.2$, thereby giving in this special case an affirmative answer to a question raised by Rauf in \cite{asia2}.

%\textbf{Aknowledgements}. The author would like to express his gratitude to Professor Dorin Popescu for his valuable suggestions during the writing of this article.

\section{Preliminary results}

We recall the following result of Herzog, Vladoiu and Zheng.

\begin{prop}(\cite[Proposition 3.4]{hvz})
Let $I\subset S$ be a monomial ideal. Then: $$\sdepth(I)\geq max\{1,n-g(I)+1\}.$$
\end{prop}

In the following, we give a similar result:

\begin{prop}
Let $I\subset S$ be a monomial ideal. Then $\sdepth(S/I)\geq n-g(I)$.
\end{prop}

\begin{proof}
In order to prove, we use a strategy similar with the Janet's algorithm, see \cite{imran}. %As in the proof of \cite[Proposition 3.4]{hvz}, 
We use induction on $n\geq 1$. If $n=1$ there is nothing to prove. Denote $m=g(I)$. If $m=1$, $I$ is principal and thus $\sdepth(S/I)=n-1$. Suppose $n>1$ and $m>1$. Let $q=deg_{x_n}(I):=\max\{j:\; x_n^j|u$ for some $u\in G(I)\}$. For all $j\leq q$, we denote $I_j$ the monomial ideal in $S'=K[x_1,\ldots,x_{n-1}]$ such that $I\cap x_n^jS'=x_n^j I_j$. Note that $g(I_j)<g(I)$ for all $j<q$ and $g(I_q)\leq g(I)$. We have
\[ S/I = S'/I_0 \oplus x_n(S'/I_1) \oplus \cdots \oplus x_n^{q-1}(S'/I_{q-1})\oplus x_n^q (S'/I_q)[x_n].\]
It follows that $\sdepth(S/I)\geq \min\{\sdepth(S'/I_j),j<q, \sdepth(S'/I_q)+1\}$. By induction hypothesis, it follows that
$\sdepth(S'/I_j)\geq n-1-q(I_j) \geq n-1-(m-1)=n-m$ for all $j<q$. Also, $\sdepth(S'/I_q)\geq n-1-g(I_q)\geq n-1-m$. This completes the proof.
\end{proof}

For any monomial ideal $J\subset S$, we denote $J^c$ the $K$-vector space spanned by all the monomials not contained in $J$. With this notation, we have the following lemma.

\begin{lema}
Let $I\subset S=K[x_1,\ldots,x_n ]$ be a monomial ideal and $v\in S$ a monomial. Then
$I = ((v)^c\cap I)\oplus v(I:v)$ and $I^c = ((v)^c\cap I^c)\oplus v(I:v)^c$.
\end{lema}

\begin{proof}
We have $I = S \cap I = ((v)^c\oplus (v)) \cap I = ((v)^c\cap I)\oplus ((v)\cap I)$. In order to complete the proof, it is enough to show that $((v)\cap I)=v(I:v)$. Indeed, if $w\in (v)\cap I$ is a monomial, then $w=vy$ for some monomial $y\in S$. Moreover, since $vy=w\in I$ it follows that $y\in (I:v)$ and thus $w\in v(I:v)$. The inclusion $((v)\cap I) \supseteq v(I:v)$ is similar. Analogously, we prove the second statement.
\end{proof} \pagebreak

\begin{teor}
Let $I\subset S$ be a monomial ideal which is not principal. Assume $I=vI'$, where $v\in S$ is a monomial and $I'=(I:v)$. Then:

(1) $\sdepth(S/I)=\sdepth(S/I')$.

(2)	$\sdepth(I)=\sdepth(I')$.	
\end{teor}

\begin{proof}
(1) By Lemma $1.3$, $S/I = I^c=(v)^c \oplus v(I'^c)$. Given a Stanley decomposition \linebreak $S/I'=\bigoplus_{i=1}^r u'_iK[Z_i]$ of $S/I'$, it follows that $\bigoplus_{i=1}^r vu'_iK[Z_i]$ is a Stanley decomposition of $v(I'^c)$. On the other hand, since $\sdepth(S/(v))=1$, one can easily give a Stanley decomposition $\mathcal D$ of $(v)^c$ with $\sdepth(\mathcal D)=n-1$. Thus, we obtain a Stanley decomposition of $S/I$ with its Stanley depth $\geq \sdepth(S/I')$. It follows that $\sdepth(S/I)\geq \sdepth(S/I')$.

In order to prove the converse inequality, we consider $S/I=\bigoplus_{i=1}^r u_i K[Z_i]$ a Stanley decomposition of $S/I$. It follows that $v(I'^c) = \bigoplus_{i=1}^r (u_i K[Z_i]\cap v(I'^c)) = \bigoplus_{i=1}^r (u_i K[Z_i]\cap (v))$. The last equality follows from the fact that $I^c\cap (v) = v(I'^c)$. We claim that $u_i K[Z_i]\cap (v)\neq (0)$ implies $LCM(u_i,v)\in u_iK[Z_i]$. Indeed, if $LCM(u_i,v)\notin u_i K[Z_i]$ it follows that $v/GCD(u_i,v)\notin K[Z_i]$ and therefore, there exists $x_j| v/GCD(u_i,v)$ such that $x_j\notin Z_i$. Thus, $v$ cannot divide any monomial of the form $u_iy$, where $y\in K[Z_i]$ and therefore $u_i K[Z_i]\cap (v)=(0)$, a contradiction. Now, since $LCM(u_i,v)\in u_i K[Z_i]$, it follows that $LCM(u_i,v)K[Z_i]\subset u_iK[Z_i]$. Obviously, $LCM(u_i,v)K[Z_i]\subset (v)$ and thus $LCM(u_i,v)K[Z_i]\subset u_i K[Z_i]\cap (v)$. On the other hand, if $u\in u_i K[Z_i]\cap (v)$ is a monomial, it follows that $u_i|u$ and $v|u$ and therefore, $LCM(u_i,v)|u$. Since $u\in u_i K[Z_i]$, it follows that $u=u_i w_i$, where  $supp(w_i)\subset Z_i$. Moreover, $supp(u/LCM(u_i,v))\subset Z_i$ and thus, $u\in LCM(u_i,v)K[Z_i]$. We obtain $LCM(u_i,v)K[Z_i] = u_i K[Z_i]\cap (v)$. In conclusion, \[ vI'^c = \bigoplus_{(v)\cap u_iK[Z_i]\neq 0}LCM(u_i,v)K[Z_i]\;so\; I'^c = \bigoplus_{(v)\cap u_iK[Z_i]\neq 0}\frac{u_i}{GCD(u_i,v)}K[Z_i]. \] It follows that $\sdepth(S/I')\geq \sdepth(S/I)$, as required.

(2) Follows from the linear space isomorphism $I'\cong vI' = I$.
%Suppose $I=\bigoplus_{i=1}^r u_iK[Z_i]$ is a Stanley decomposition for $I$. One can easily see that $I'=\bigoplus_{i=1}^r u_i/v K[Z_i]$ is a Stanley decomposition for $I'$, and thus $\sdepth(I)\leq \sdepth(I')$. Conversely, if $I'=\bigoplus_{i=1}^r \bar{u}_i K[Z_i]$ is a Stanley decomposition of $I'$ it follows that \linebreak $I=\bigoplus_{i=1}^r \bar{u}_iv K[Z_i]$ is a Stanley decomposition of $I$ and thus $\sdepth(I')\leq \sdepth(I)$.
\end{proof}

\begin{prop}
Let $I\subset S$ be a monomial ideal. Then:

(1)	$\sdepth(S/I)\geq n-c(I)$

(2)	$\sdepth(I)\geq n-c(I)+1$.
\end{prop}

\begin{proof}
(1) Let $v=GCD(u|\;u\in G(I))$ and $I'=(I:v)$. By $1.4(1)$ we can assume that $I'=I$. By reordering the variables, we can assume that $I\subset (x_1,x_2,\ldots,x_m)$, where $m=c(I)$. We write $I=(I\cap K[x_1,\ldots,x_m])S$. \cite[Lemma 3.6]{hvz} implies $\sdepth(S/I) = \sdepth(K[x_1,\ldots,x_m]/(I\cap K[x_1,\ldots,x_m])) + n-m \geq n-m$.
(2) The proof is similar, if we see that $\sdepth(I\cap K[x_1,\ldots,x_m])\geq 1$, see $1.1$.
\end{proof}

\begin{prop}
Let $I\subset S$ be a monomial ideal which is not principal with $c(I)=2$ or $g(I)=2$. Then $\sdepth(I)=n-1$ and $\sdepth(S/I)=n-2$.
\end{prop}

\begin{proof}
If $c(I)=2$, then, by $1.5(2)$, it follows that $\sdepth(I)\geq n-c(I)+1=n-1$. If $g(I)=2$, by $1.1$, it follows that $\sdepth(I)\geq n-1$. On the other hand, $\sdepth(I)<n$, otherwise, $I$ would be principal. Thus $\sdepth(I)=n-1$.

According to $1.5(1)$ or $1.2$, $\sdepth(S/I)\geq n-2$ if $c(I)=2$ or, respectively, $g(I)=2$. We consider the case $c(I)=2$. Let $v=GCD(u|\;u\in G(I))$ and $I'=(I:v)$. By $1.4(1)$, we can assume that $I=I'$ and $supp(I)=\{x_1,x_2\}$. Since $v=1$, it follows that $x_1^a, x_2^b\in G(I)$ for some positive integers $a$ and $b$. Therefore, $\sdepth(I\cap K[x_1,x_2])=0$ and moreover,\linebreak \cite[Lemma 3.6]{hvz} implies $\sdepth(S/I) = \sdepth(K[x_1,x_2]/(I\cap K[x_1,x_2])) + n-2 = n-2$.

%Let $w=x_1^{a-1}$. Obviously, $w\notin I$, but $x_1w\in I$ and $x_2^kw\in I$ for $k\geq b$. It follows that $w$ cannot be contained in a Staley space of $S/I$ with dimension $\geq n-1$. Thus $\sdepth(S/I)=n-2$.

We consider now the case $g(I)=2$. Suppose $I=(u_1,u_2)$. By $1.4(1)$, we can assume $GCD(u_1,u_2)=1$. Therefore, $I$ is a complete intersection and by \cite[Proposition 1.2]{hsy} or \cite[Corollary 1.4]{asia}, it follows that $\sdepth(S/I)=n-2$.
\end{proof}

\section{Stanley depth of monomial ideals with small number of generators}

The following result is a particular case of \cite[Theorem 1.4]{cim}.

\begin{teor}
Let $I\subset S$ be a monomial ideal. Then $\sdepth(S/I)=0$ if and only if $depth(S/I)=0$.
\end{teor}

Note that $depth(S/I)=0$ if and only if $I\neq I^{sat}$, where $I^{sat}=\bigcup_{k\geq 1}(I:(x_1,\ldots,x_n)^k)$ is the \emph{saturation} of $I$.

\begin{cor}
Let $I\subset S$ be a monomial ideal. Then $\sdepth(S/I)=0$ if and only if $\sdepth(S/I^k)=0$, where $k\geq 1$.
\end{cor}

\begin{proof}
It is enough to notice that $\depth(S/I)=0$ if and only if $\depth(S/I^k)=0$, where $k\geq 1$, than we apply Theorem $1.1$.
\end{proof}

\begin{cor}
Let $I\subset S$ be a monomial ideal with $c(I)=n$ and $(x_1,\ldots,x_{n-1})\subset \sqrt{I}$. Then $\sdepth(S/I)=0$.
\end{cor}

\begin{proof}
Since $(x_1,\ldots,x_{n-1})\subset \sqrt{I}$ it follows that for all $1\leq j\leq n-1$, there exists a positive integer $a_j$ such that $x_j^{a_j}\in I$. Since $c(I)=n$ it follows that there exists a monomial $u\in G(I)$ with $x_n|u$. If $u=x_n^{a_n}$, it follows that $I$ is artinian and thus, by Theorem $2.1$, $\sdepth(S/I)=0$. Suppose this is not the case. We consider $w=u/x_n$. Obviously, $x_j^{a_j}w\in I$ for any $1\leq j\leq n$, where $a_n:=1$. Thus, $w\in I^{sat}\setminus I$ and then $\sdepth(S/I)=0$ by $2.1$.
\end{proof}

\begin{teor}
Let $I\subset S$ be a monomial ideal with $g(I)=3$. Then $\sdepth(I)=n-1$.
\end{teor}

\begin{proof}
Denote $G(I)=\{v_1,v_2,v_3\}$. By $1.4(2)$, we can assume that $GCD(v_1,v_2,v_3)=1$. If $I$ is a complete intersection, by \cite[Proposition 3.8]{hvz} or \cite[Theorem 2.4]{shen}, it follows that $\sdepth(I)=n-1$. If this is not the case, it follows that there exists a variable, let's say $x_n$, such that $x_n|v_1$,$x_n|v_2$ and $x_n$ does not divide $v_3$. Let $a=deg_{x_n}(v_1)$ and $b=deg_{x_n}(v_2)$ and suppose $a\leq b$. We have the following decomposition given by Lemma $1.3$: 
\[ I = (I\cap (x_n^a)^c) \oplus x_n^a (I:x_n^a) =  \bigoplus_{j=0}^{a-1} x_n^j v_3 K[x_1,\ldots,x_{n-1}] \oplus x_n^a (I:x_n^a). \]
Note that $g((I:x_n^a))\leq g(I)$. If $g((I:x_n^a))<3$, we can find a Stanley decomposition for $(I:x_n^a)$ with it's Stanley depth $\geq n-1$ and we stop. Otherwise, we replace $I$ with $(I:x_n^a)$ and we repeat the previous procedure until we obtain an ideal with $\leq 2$ generators of a monomial complete intersection ideal. Finally, we obtain a Stanley decomposition of $I$ with it's Stanley depth equal to $n-1$. On the other hand, $\sdepth(I)<n$, since $I$ is not principal.
\end{proof}

\begin{exm}
Let $I=(x_1^3,x_2^2x_3^2,x_1x_2^3x_3)$. We have:
\[ I = ((x_2^{2})^{c}\cap I)\oplus x_2^{2}(I:x_2^{2})=x_1^{3} K[x_1,x_3]\oplus x_1^{3}x_2 K[x_1,x_3]\oplus x_2^2(x_1^3,x_3^2,x_1x_2x_3),\]
On the other hand,
\[ (x_1^3,x_3^2,x_1x_2x_3) = x_3^2 K[x_2,x_3] \oplus x_1 (x_1^2,x_3^2,x_2x_3) = x_3^2 K[x_2,x_3] \oplus x_1^3K[x_1,x_2]\oplus x_1x_3 (x_1^2,x_2,x_3).\] We obtain the following Stanley decomposition of $I$:
\[ I = x_1^{3} K[x_1,x_3]\oplus x_1^{3}x_2 K[x_1,x_3]\oplus x_2^2x_3^2 K[x_2,x_3] \oplus x_2^2x_1^3K[x_1,x_2] \oplus \]
\[ \oplus x_1x_2^2x_3(x_2K[x_2,x_3]\oplus x_1x_2K[x_2,x_3] \oplus x_1^2K[x_1,x_2]\oplus x_3K[x_1,x_3]\oplus x_1^2x_2x_3K[x_1,x_2,x_3] ).\]
\end{exm}

\begin{teor}
Let $I\subset S$ be a monomial ideal with $g(I)\leq 3$. Then, $I$ and $S/I$ satisfy the Stanley's conjecture.
\end{teor}

\begin{proof}
It is well known that $\depth(I)=\depth(S/I)+1$ and $\depth(S/I)=n-1$ if and only if $I$ is principal. Thus, if $g(I)=1$, there is nothing to prove. If $g(I)=2$, by $1.6$, we have $\sdepth(S/I)=n-2$ and $\sdepth(I)=n-1$. Since $I$ is not principal, it follows that $\depth(S/I)\leq n-2$ and $\depth(I)\leq n-1$. So we are done.

We consider the case $g(I)=3$. By $2.4$, $\sdepth(I)=n-1$ and thus $\sdepth(I)\geq \depth(I)$. According to $1.2$, $\sdepth(S/I)\geq n-3$. Thus, if $\depth(S/I)\leq n-3$ we are done. 

Now, assume $\depth(S/I)=n-2$. Denote $G(I)=\{v_1,v_2,v_3\}$, $v=GCD(v_1,v_2,v_3)$ and $I'=(I:v)$. By $1.4(1)$, $\sdepth(S/I) = \sdepth(S/I')$. On the other hand, by \cite[Corollary 1.3]{asia2}, $\depth(S/I'))\geq \depth(S/I)$. In fact, $\depth(S/I')=n-2$, since $I'$ is not principal. Thus, we can assume $I = I'$. Note that $I$ is not a complete intersection, otherwise $\depth(S/I)=n-3$. Therefore, there exists a variable, let's say $x_n$, such that $x_n|v_1$,$x_n|v_2$ and $x_n$ does not divide $v_3$. Let $a=deg_{x_n}(v_1)$ and $b=deg_{x_n}(v_2)$ and suppose $a\leq b$. We have the following decomposition given by Lemma $1.3$: 
\[ S/I = I^c = (I^c \cap (x_n^a)^c) \oplus x_n^a ((I:x_n^a)^c) =  \bigoplus_{j=0}^{a-1} x_n^j (S'/ v_3S') \oplus x_n^a ((I:x_n^a)^c), \]
where $S'=K[x_1,\ldots,x_{n-1}]$. Note that $\sdepth_{S'}{S'/v_3S'}=(n-1)-1=n-2$. We have $g((I:x_n^a))\leq g(I)$. If $g((I:x_n^a))<3$, we can find a Stanley decomposition for $S/(I:x_n^a)$ with it's Stanley depth $>=n-2$ and we stop.
Assume $g((I:x_n^a))=3$. By \cite[Corollary 1.3]{asia2}, we have $\depth(S/(I:x_n^a)) \geq \depth(S/I) = n-2$. Thus, $\depth(S/(I:x_n^a))=n-2$ and therefore $(I:x_n^a)$ is not a complete intersection. We replace $I$ with $(I:x_n^a)$ and we continue the previous procedure. Finally, we obtain a Stanley decomposition for $S/I$ with it's Stanley depth equal to $n-2$. Therefore $\sdepth(S/I)\geq n-2$, as required.
\end{proof}

\section{Monomial ideals in three variables}

\begin{lema}
Let $I\subset S:=K[x_1,x_2,x_3]$ be a monomial ideal with $v:=GCD(u|u\in G(I))=1$. For $1\leq j\leq 3$, we denote $S_j:=K[Z_j]$ and $I_j=I\cap S_j$, where $Z_j=\{x_1,x_2,x_3\}\setminus \{x_j\}$. If $I^{sat} = I$ then there exists some $1\leq j\leq 3$ such that $I_j^{sat}=I_j$. 
\end{lema}

\begin{proof}
If $I=S$ there is nothing to prove, so we can assume $I\neq S$. Since $I^{sat}=I$, it follows that $\mathbf{m}=(x_1,x_2,x_3)\notin Ass(S/I)$. Since $v=1$, it follows that $(x_j)\notin Ass(S/I)$ for all $1\leq j\leq 3$. We denote $\mathbf{m}_j$ the ideal generated by the variables from $Z_j$. We have $Ass(S/I)\subset \{\mathbf{m}_1,\mathbf{m}_2,\mathbf{m}_3\}$. 

Thus, we can find a decomposition $I=\bigcap_{j=1}^3 Q_j$ such that $Q_j$ is $\mathbf{m}_j$-primary or $Q_j=S$ for all $1\leq j\leq 3$. It follows that $I_k=\bigcap_{j=1}^3 (Q_j\cap S_k)$. We assume $Q_1=(x_2^a,x_3^b,\ldots)$, $Q_2=(x_1^c,x_3^d,\ldots)$ and $Q_3=(x_1^e,x_2^f,\ldots)$ where $a,b,c,d,e,f$ are some nonnegative integers. By reordering $Q_i$'s, we can assume $a\geq f$. It follows that $Q_1\cap S_3 = (x_2^a) \subset Q_3$. Therefore, $I_3 = (x_2^a)\cap (x_1^c) = (x_2^ax_1^c)$ is principal and thus, $I_3=I_3^{sat}$.
\end{proof}

\begin{prop}
Let $I\subset S:=K[x_1,x_2,x_3]$ be a monomial ideal which is not principal. If $I=I^{sat}$ then $\sdepth(I)=2$.
\end{prop}

\begin{proof}
We denote $v=GCD(u|u\in G(I))$ and $I'=(I:v)$. By Theorem $1.4(2)$, we have $\sdepth(I)=\sdepth(I')$. Since, also, $I^{sat} = vI'^{sat}$, we can assume $I=I'$. If $c(I)=2$ or $g(I)=2$, by Proposition $1.6$, it follows that $\sdepth(I)=2$. 
Now, we consider the case: $c(I)=3$ and $g(I)\geq 3$. In the notations of Lemma $3.1$, we can assume that $I_1^{sat}=I_1$ and $I_1$ is principal. Thus $\sdepth(I_1)=2$. We write
$I=I_1\oplus x_1(I:x_1)$. Obviously, $I\subsetneq (I:x_1)$ and $(I:x_1)^{sat}=(I:x_1)$. We can use the same procedure for $(I:x_1)$. Finally, we obtain a Stanley decomposition of $I$ with its Stanley depth equal to $2$.
\end{proof}

\begin{cor}
If $I\subset K[x_1,x_2,x_3]$ is a monomial ideal, then $\sdepth(I) \geq \sdepth(S/I)+1$. In particular, if $\sdepth(I)=1$, then $\depth(I)=1$.
\end{cor}

\begin{proof}
If $\sdepth(S/I)=0$ there is nothing to prove, since $\sdepth(I)\geq 1$. If $\sdepth(S/I)=1$, by Theorem $2.1$, it follows that $I=I^{sat}$ and, by $3.1$, $\sdepth(I)=2$. On the other hand, $\sdepth(S/I)=2$ if and only if $I$ is principal and, thus, if and only if $\sdepth(I)=3$.

For the second statement, assume $\depth(I)>1$, i.e. $\depth(S/I)>0$. It follows by $2.1$ that $I=I^{sat}$ and thus, by $3.2$, $\sdepth(I)\geq 2$, a contradiction.
\end{proof}

\begin{obs}
\emph{A similar result to Lemma $3.1$ is not true for $n\geq 4$. Let $S=K[x_1,x_2,x_3,x_4]$, $Q_1=(x_2^3,x_3^2,x_4)$, $Q_2=(x_1^3,x_3,x_4^2)$, $Q_3=(x_1^2,x_2,x_4^3)$, $Q_4=(x_1,x_2^2,x_3^3)$ and $I=Q_1\cap Q_2\cap Q_3\cap Q_4$. Let $Z_k=\{x_1,x_2,x_3,x_4\}\setminus \{x_k\}$ and $S_k:=K[Z_k]$, where $1\leq k\leq 4$. One can easily see that $I_k=I\cap S_k = \bigcap_{j=1}^4 (Q_j\cap S_k)$ is a reduced primary decomposition of $I_k$. In particular, $\mathbf{m}_k = \sqrt{Q_k\cap S_k} \in Ass(S_k/I_k)$ and thus $I_k^{sat}\neq I_k$. On the other hand, $I=I^{sat}$.}
\end{obs}

\vspace{2mm} \noindent {\footnotesize
\begin{minipage}[b]{15cm}
 Mircea Cimpoeas, Institute of Mathematics of the Romanian Academy, Bucharest, Romania\\
 E-mail: mircea.cimpoeas@imar.ro
\end{minipage}}

\begin{thebibliography}{99}%{breitestes Label}  
  \bibitem[1]{pops}Sarfraz Ahmad, Dorin Popescu "Sequentially Cohen-Macaulay monomial ideals of embedding dimension four", Bull. Math. Soc. Sc. Math. Roumanie 50(98), no.2 (2007), p.99-110.
  \bibitem[2]{imran}Imran Anwar "Janet's algorithm", Bull. Math.
                    Soc. Sc. Math. Roumanie 51(99), no.1 (2008), p.11-19.
  \bibitem[3]{popi}Imran Anwar, Dorin Popescu "Stanley Conjecture in small embedding dimension", Journal of Algebra 318 (2007),                            p.1027-1031.
  \bibitem[4]{apel}J.Apel "On a conjecture of R.P.Stanley", Journal of Algebraic Combinatorics, 17(2003), p.36-59.
  %\bibitem[5]{par}Csaba Biro, David M.Howard, Mitchel T.Keller, William T.Trotter, Stephen J.Young "Partitioning subset lattices into intervals, preliminary version", Preprint 2008.   
  \bibitem[5]{mir}Mircea Cimpoeas "Stanley depth for monomial complete intersection", Bull. Math.
                  Soc. Sc. Math. Roumanie 51(99), no.3 (2008), p.205-211.
  \bibitem[6]{cim}Mircea Cimpoeas "Some remarks on the Stanley depth for multigraded modules", Le Mathematiche, Vol. LXIII (2008) – Fasc. II, pp. 165–171. 
  \bibitem[7]{hsy}J\"urgen Herzog, Ali Soleyman Jahan, Siamak Yassemi "Stanley decompositions and partitionable simplicial complexes",                    Journal of Algebraic Combinatorics 27(2008), p.113-125.
  \bibitem[8]{hvz}J\"urgen Herzog, Marius Vladoiu, Xinxian Zheng "How to compute the Stanley depth of a monomial ideal", Journal of Algebra 2009, doi:10:1016/j.jalgebra.2008.01.006 in press %, http://arxiv.org/pdf/0712.2308.   
  \bibitem[9]{jah}Ali Soleyman Jahan "Prime filtrations of monomial ideals and polarizations", Journal of Algebra  312 (2007),                            p.1011-1032.
  \bibitem[10]{sum}Sumiya Nasir "Stanley decompositions and localization", 
                  Bull. Math. Soc. Sc. Math. Roumanie 51(99), no.2 (2008), p.151-158.
  \bibitem[11]{pope}Dorin Popescu "Stanley depth of multigraded modules", Journal of Algebra, 321 (10), 2009, p.2782-2797.
  \bibitem[12]{asia}Asia Rauf "Stanley decompositions, pretty clean filtrations and reductions modulo regular
   elements", Bull. Soc. Sc. Math. Roumanie 50(98), no.4 (2007), p.347-354.
  \bibitem[13]{asia2}Asia Rauf "Depth and Stanley depth of multigraded modules", Preprint 2009 %http://arxiv.org/pdf/0812.2080
  , to appear in Communications in Algebra.
  \bibitem[14]{shen}Yihuang Shen "Stanley depth of complete intersection monomial ideals and upper-discrete
    partitions", Journal of Algebra 321(2009), 1285-1292. 
\end{thebibliography}
\end{document}